\documentclass{amsart}

\usepackage{amssymb}

\newtheorem{theorem}{Theorem}[section]
\newtheorem{lemma}[theorem]{Lemma}
\newtheorem{proposition}[theorem]{Proposition}

\theoremstyle{definition}
\newtheorem{definition}[theorem]{Definition}

\theoremstyle{remark}
\newtheorem{remark}[theorem]{Remark}

\newcommand{\N}{\mathbb{N}}
\newcommand{\Z}{\mathbb{Z}}
\newcommand{\R}{\mathbb{R}}
\newcommand{\C}{\mathbb{C}}
\newcommand{\T}{\mathbb{T}}

\newcommand{\G}{\widehat{G}}

\newcommand{\g}{\mathfrak{g}}
\newcommand{\h}{\mathfrak{h}}
\newcommand{\mt}{\mathfrak{T}}

\newcommand{\Fourier}{\mathcal{F}}
\newcommand{\Matrix}{\mathcal{M}}
\newcommand{\Lebesgue}{\mathcal{L}}
\newcommand{\Constant}{\mathcal{C}}
\newcommand{\Ball}{\mathcal{U}}
\newcommand{\Root}{\mathcal{R}}
\newcommand{\Weyl}{\mathcal{W}_G}

\newcommand{\MT}{\textnormal{\textrm{\textbf{T}}}}
\newcommand{\Identity}{\mathbf{1}}
\newcommand{\Chamber}{\mathbf{C}}
\newcommand{\Simple}{\emph{1}}
\newcommand{\Planes}{\mathrm{P}}

\begin{document}

\title[Sharp Fourier type and cotype]
{Sharp Fourier type and cotype \\ with respect to compact
semisimple Lie groups}

\author[Garc\'{\i}a-Cuerva, Marco and Parcet]
{Jos\'{e} Garc\'{\i}a-Cuerva, Jos\'{e} Manuel Marco and Javier Parcet}

\address{Department of Mathematics, Universidad Aut\'{o}noma de
Madrid, Madrid 28049, Spain}

\email{jose.garcia-cuerva@uam.es} \email{javier.parcet@uam.es}

\thanks{Research supported in part by the European Commission
via the TMR Network `Harmonic Analysis' and by Project BFM
2001/0189, Spain}


\date{}


\begin{abstract}
Sharp Fourier type and cotype of Lebesgue spaces and Schatten
classes with respect to an arbitrary compact semisimple Lie group
are investigated. In the process, a local variant of the
Hausdorff-Young inequality on such groups is given.
\end{abstract}

\maketitle

\section*{Introduction}

Let $1 \le p \le 2$. An operator space $E$ is said to have Fourier
type $p$ with respect to the compact group $G$ if the
vector-valued Fourier transform extends to a completely bounded
map $$\Fourier_{G,E}: L_E^p(G) \longrightarrow
\Lebesgue_E^{p'}(\G)$$ where $p' = p / (p-1)$ is the exponent
conjugate to $p$. That is, a vector-valued Hausdorff-Young
inequality of exponent $p$ is satisfied. Similarly, if we replace
the operator $\Fourier_{G,E}$ by its inverse, we get the notion of
Fourier cotype $p'$ of $E$ with respect to $G$. Following the
notation of \cite{GP}, we define the constants
$$\Constant_p^1(E,G) = \|\Fourier_{G,E}\|_{cb(L_E^p(G),
\Lebesgue_E^{p'}(\G))} \qquad \mbox{and} \qquad
\Constant_{p'}^2(E,G) =
\|\Fourier_{G,E}^{-1}\|_{cb(\Lebesgue_E^p(\G), L_E^{p'}(G))}.$$

The Fourier type and cotype become stronger conditions on the pair
$(E,G)$ as the exponents $p$ and $p'$ approach $2$. This gives
rise to the notions of sharp Fourier type and cotype exponents.
The present paper grew out of the project to investigate the sharp
Fourier type and cotype of Lebesgue spaces $L^p$ and Schatten
classes $S^p$ and it is a natural continuation of \cite{GP}.
However, as we shall see below, some other results have appeared
in the process which are interesting on their own right.

In section \ref{Section-Statement} we recall that the natural
candidates for the sharp Fourier type and cotype of $L^p$ and
$S^p$ --where now $1 \le p \le \infty$-- are $\min(p,p')$ and
$\max(p,p')$ res\-pectively. To justify that this guess is right,
one would have to show that for $1 \le p < q \le 2$
$$\begin{array}{llclclc} (a) & \Constant_q^1(L^p(\Omega),G) & = &
\Constant_{q'}^2(L^{p'}(\Omega),G) & = & \infty & \\ (b) &
\Constant_q^1(L^{p'}(\Omega),G) & = &
\Constant_{q'}^2(L^p(\Omega),G) & = & \infty & \qquad \qquad
\qquad \qquad \qquad \qquad \end{array}$$ with the obvious
modifications for the Schatten classes. Section
\ref{Section-Remarks} is devoted to make some remarks about $(a)$
and $(b)$. First we show that, to have any chance of getting
positive answers to these questions, we have to require the group
$G$ not to be finite and the operator spaces $L^p$ and $S^p$ to be
infinite-dimensional. Then, under such assumptions, one can easily
get the following inequality $$\Constant_q^1(L^p(\Omega),G) \ge
\limsup_{n \rightarrow \infty} \,\ \Constant_q^1(l^p(n),G)$$ and
the analog for $L^{p'}(\Omega)$. Therefore the growth of
$\Constant_q^1(l^p(n),G)$ and $\Constant_q^1(l^{p'}(n),G)$
provides a possible way to obtain $(a)$ and $(b)$. In the last
part of section \ref{Section-Remarks} we analyze the vector-valued
Lebesgue spaces and Schatten classes.

The growth of $\Constant_q^1(l^p(n),G)$ is investigated in section
\ref{Section-p<2}. To be precise, if $G$ stands for a compact
semisimple Lie group and $1 \le p < q \le 2$, then there exists a
constant $0 \le \mathcal{K}(G,q) \le 1$ such that
$\Constant_q^1(l^p(n),G) \ge \mathcal{K}(G,q) \,\ n^{1/p - 1/q}$
for all $n \ge 1$. If one is able to see that $\mathcal{K}(G,q) >
0$, this result gives $(a)$. Moreover, we would obtain optimal
growth since $\Constant_q^1(l^p(n),G) \le n^{1/p - 1/q}$ for any
compact group. We shall see that $$\mathcal{K}(G,q) = \inf_{n \ge
1} \sup \left\{
\frac{\|\widehat{f}\|_{\Lebesgue^{q'}(\G)}}{\|f\|_{L^q(G)}}: \,\ f
\ \ \mbox{central}, \ \ f \in L^q(G), \ \ \mbox{supp}(f) \subset
\Ball_n \right\}$$ where $\{\Ball_n: n \ge 1\}$ is a basis of
neighborhoods of $\Identity$, the identity element of $G$. The
Hausdorff-Young inequality on compact groups provides
$\mathcal{K}(G,q) \le 1$. The interesting point lies in the
inequality $\mathcal{K}(G,q) > 0$ which constitutes a local
variant of the Hausdorff-Young inequality on $G$ with exponent
$q$.

Sections \ref{Section-Simple} and \ref{Section-Local} are
completely devoted to the proof of this local inequality. In the
abelian setting, the particular case $G = \T$ was explored by
Andersson in \cite{An}. The basic idea is to consider a function
$f: \T \rightarrow \C$ as a complex-valued function on $\R$
supported in $[-1/2, 1/2)$. Then, by expressing the norm of
$\widehat{f}$ on $L^{q'}(\R)$ as a Riemann sum, one obtains
$$\frac{\|\widehat{f}\|_{L^{q'}(\R)}}{\|f\|_{L^q(\R)}} = \lim_{k
\rightarrow \infty}
\frac{\|\widehat{\varphi}_k\|_{L^{q'}(\T)}}{\|\varphi_k\|_{L^q(\T)}}$$
where $\varphi_k(t) = k^{1/q} f(kt)$. This gives
$\mathcal{K}(\T,q) \ge \mathcal{B}_q$ --where $\mathcal{B}_q =
\sqrt{q^{1/q}/q^{'1/q'}}$ stands for the constant of
Babenko-Beckner, see \cite{Ba} and \cite{B1}-- but in fact the
equality holds, as it was proved by Sj\"{o}lin in \cite{Sj}. We show
here that Andersson's argument, suitably modified, is also valid
in the context of compact semisimple Lie groups. In section
\ref{Section-Simple} we summarize the main results of the
structure and representation theory of compact semisimple Lie
groups that will be used in the process. Then we use these
algebraic results to get an expression for the Fourier transform
of central functions $f: G \rightarrow \C$ in terms of the Fourier
transform $\Fourier_{\MT}$ on the maximal torus $\MT$ of $G$. This
will allow us to work over the maximal torus where we know that
Andersson obtained a satisfactory result. However, in the
non-commutative setting, the degree $d_{\pi}$ of an irreducible
representation $\pi$ does not have to be $1$. And we shall see
that this becomes a further obstacle to be treated in section
\ref{Section-Local}. There we combine some results --as the Weyl
dimension formula-- concerning the representation theory of
compact semisimple Lie groups with classical harmonic analysis to
avoid this difficulty.

On the other hand, if we notice that $\Constant_q^1(l^{p'}(n),G) =
\Constant_{q'}^2(l^p(n),G)$, we can understand the growth of this
constant as the dual problem of the growth of
$\Constant_q^1(l^p(n),G)$ in the sense that we substitute the
Fourier transform operator $\Fourier_{G, l^p(n)}$ by its inverse.
Therefore, since the dual object is no longer a group --as it is
when $G$ is abelian-- we do not have a Fourier inversion theorem
and we should not expect to reconstruct the proof given in
sections \ref{Section-p<2}, \ref{Section-Simple} and
\ref{Section-Local} step by step. At the time of this writing, we
are not able to solve this problem and so we pose it as follows:

\

\noindent \textbf{Problem:} Let $G$ be any compact semisimple Lie
group and let $1 \le p < q \le 2$. Does the estimate
$\Constant_q^1(l^{p'}(n),G) \ge \mathcal{K}(G,q) \,\ n^{1/p -
1/q}$ hold for some positive constant $\mathcal{K}(G,q)$ depending
only on $G$ and $q$?

\

Finally, we point to a non-commutative notion of Rademacher type
for operator spaces, see \cite{GP2}. We think this notion could be
helpful in order to study the growth of
$\Constant_q^1(l^{p'}(n),G)$.

\section{Statement of the problem}
\label{Section-Statement}

All throughout this paper some basic notions of operator space
theory and non-commutative vector-valued discrete $L^p$ spaces
will be assumed. The definitions and results about operator spaces
that we are using can be found in the book of Effros and Ruan
\cite{ER2}, while for the study of our non-commutative $L^p$
spaces the reader is referred to \cite{P2}, where Pisier analyzes
them in detail. In any case all the analytic preliminaries of this
paper are summarized in \cite{GP}, where we study the Fourier type
and cotype of an operator space with respect to a compact group.
In order to state the problem we want to solve, we begin by
recalling the definitions and the main properties of Fourier type
and cotype.

Let $G$ be a compact topological group endowed with its Haar
measure $\mu$ normalized so that $\mu(G) = 1$ and let $\pi \in \G$
be an irreducible unitary representation of $G$ of degree
$d_{\pi}$. Here the symbol $\G$ stands for the dual object of $G$.
Given an operator space $E$, it was shown in \cite{GP} that --by
fixing a basis on the representation space of each $\pi \in \G$--
the Fourier transform operator $\Fourier_{G,E}$ for functions
defined on $G$ and with values on $E$, has the form $f \in
L_E^1(G) \longmapsto \big( \widehat{f}(\pi) \big)_{\pi \in \G} \in
\Matrix_E(\G)$, where $$\widehat{f}(\pi) = \int_G f(g)
\pi(g)^{\star} d \mu(g) \qquad \qquad \mbox{and} \qquad \qquad
\Matrix_E(\G) = \prod_{\pi \in \G} M_{d_{\pi}} \otimes E.$$ Here
$M_n$ denotes the space of $n \times n$ complex matrices. Let $1
\le p < \infty$, if $S_n^p(E)$ stands for the vector-valued
Schatten class on $M_n \otimes E$, we define the spaces
\begin{itemize}
\item $\displaystyle \Lebesgue_E^p(\G) = \Big\{A \in
\Matrix_E(\G): \,\ \|A\|_{\Lebesgue_E^p(\G)} = \Big( \sum_{\pi \in
\G} d_{\pi} \|A^{\pi}\|_{S_{d_{\pi}}^p(E)}^p \Big)^{1/p} < \infty
\Big\}$
\item $\displaystyle \Lebesgue_E^{\infty}(\G) = \Big\{A \in
\Matrix_E(\G): \,\ \|A\|_{\Lebesgue_E^{\infty}(\G)} = \sup_{\pi
\in \G} \|A^{\pi}\|_{S_{d_{\pi}}^{\infty}(E)} < \infty \Big\}$.
\end{itemize}
We write $\Lebesgue^p(\G)$ for the case $E = \C$. Finally, let $1
\le p \le 2$, then by the Hausdorff-Young inequality on compact
groups --see \cite{GP} or Kunze's paper \cite{Ku}-- it is not
difficult to check that $\mathcal{F}_{G,E} (L^p(G) \otimes E)
\subset \Lebesgue^{p'}(\G) \otimes E$ and $\mathcal{F}_{G,E}^{-1}
(\Lebesgue^p(\G) \otimes E) \subset L^{p'}(G) \otimes E$. This
motivates the following definitions.
\begin{definition} Let $1 \le p \le 2$ and let $p'$ denote
its conjugate exponent. The operator space $E$ has \emph{Fourier
type} $p$ with respect to the compact group $G$ if the Fourier
transform operator $$\mathcal{F}_{G,E}: L^p(G) \otimes E
\rightarrow \Lebesgue^{p'}(\G) \otimes E$$ can be extended to a
completely bounded operator from $L_E^p(G)$ into
$\Lebesgue_E^{p'}(\G)$. In that case $\Constant_p^1(E,G)$ will
stand for its $cb$ norm.
\end{definition}
\begin{definition} In the same fashion, the operator
space $E$ has \emph{Fourier cotype} $p'$ with respect to the
compact group $G$ if the inverse $$\mathcal{F}_{G,E}^{-1}:
\Lebesgue^p(\G) \otimes E \rightarrow L^{p'}(G) \otimes E$$ can be
extended to a completely bounded operator from $\Lebesgue_E^p(\G)$
to $L_E^{p'}(G)$. As before we shall denote its $cb$ norm by
$\Constant_{p'}^2(E,G)$.
\end{definition}

One of the properties proved in \cite{GP} is that every operator
space has Fourier type $1$ and Fourier cotype $\infty$ with
respect to any compact group. In particular, the complex
interpolation method for operator spaces --see Pisier's work
\cite{P1}-- provides the following result.

\begin{lemma} \label{Order} Let $1 \le p_1 \le p_2 \le 2$ and
assume that $E$ has Fourier type $p_2$ with respect to $G$, then
$E$ has Fourier type $p_1$ with respect to $G$. Similarly, Fourier
cotype $p_2'$ of $E$ with respect to $G$ implies Fourier cotype
$p_1'$ of $E$ with respect to $G$.
\end{lemma}

Therefore, the Fourier type and cotype become stronger conditions
on the pair $(E,G)$ as the exponent $p$ --and consequently its
conjugate $p'$-- tends to $2$. So lemma \ref{Order} gives rise to
the following definition.

\begin{definition} The \emph{sharp Fourier type and cotype
exponents} of an operator space $E$ with respect to the compact
group $G$ are defined respectively by
\begin{eqnarray*}
p_1(E,G) & = & \sup \{p \le 2: E \ \ \mbox{has Fourier type} \ \ p
\ \ \mbox{with respect to} \ \ G\} \\ p_2(E,G) & = & \inf \{p' \ge
2: E \ \ \mbox{has Fourier cotype} \ \ p' \ \ \mbox{with respect
to} \ \ G\}.
\end{eqnarray*}
If $E$ has Fourier type $p_1(E,G)$ with respect to $G$ we say that
$E$ has \emph{sharp Fourier type} $p_1(E,G)$. The \emph{sharp
Fourier cotype} of $E$ is defined analogously.
\end{definition}

In order to simplify the statement of the problem we shall need
the following lemma --see \cite{GP}-- which analyzes the Fourier
type and cotype of the dual $E^{\star}$ of an operator space $E$
with respect to a compact group $G$.

\begin{lemma} \label{Duality} Let $1 \le p \le 2$ and let $p'$ be
the conjugate exponent of $p$. Then we have the equalities
$\Constant_p^1(E^{\star},G) = \Constant_{p'}^2(E,G)$ and
$\Constant_{p'}^2(E^{\star},G) = \Constant_p^1(E,G)$.
\end{lemma}

The problem we want to investigate in this paper is how to find
out the sharp Fourier type and cotype of Lebesgue spaces and
Schatten classes. Concerning these topics we present here a result
given in \cite{GP} from which we start out. In what follows
$(\Omega, \mathcal{A}, \nu)$ will denote a $\sigma$-finite or
regular measure space and $S_{\N}^p$ the classical Schatten class
over the space of compact operators on $l^2$.

\begin{theorem} \label{Lebesgue-Schatten} Let $1 \le p \le
\infty$, then the spaces $L^p(\Omega)$, $S_n^p$ and $S_{\N}^p$
have Fourier type $\min(p,p')$ and Fourier cotype $\max(p,p')$. In
fact, the vector-valued Fourier transform --or its inverse-- is a
complete contraction in each of the cases considered.
\end{theorem}

Therefore, if we consider two exponents $p$ and $q$ such that $1
\le p < q \le 2$, we would like to find out conditions on $G$ and
$\Omega$ under which $$\begin{array}{llclclc} (a) &
\Constant_q^1(L^p(\Omega),G) & = &
\Constant_{q'}^2(L^{p'}(\Omega),G) & = & \infty & \\ (b) &
\Constant_q^1(L^{p'}(\Omega),G) & = &
\Constant_{q'}^2(L^p(\Omega),G) & = & \infty & \qquad \qquad
\qquad \qquad \qquad \qquad \end{array}$$ with the obvious
modifications for the Schatten classes.

\section{Some remarks about the problem}
\label{Section-Remarks}

In this section we shall point out some remarks about the problem
we have just stated. We begin by showing some necessary conditions
that should hold to obtain a positive answer to our question.
Second we wonder about sufficient conditions that we shall work
with along the rest of this paper. Finally we study what happens
if we consider vector-valued $L^p$ spaces --or Bochner-Lebesgue
spaces-- and vector-valued Schatten classes.

\subsection{Necessary conditions}

The first necessary condition we are talking about is on the
compact group $G$. We have to exclude \textbf{finite groups} from
our treatment since, as we shall see immediately, every operator
space $E$ has sharp Fourier type and cotype $2$ with respect to
any finite group. Anyway the next result is a bit more accurate.

\begin{proposition} \label{Finite} Let $G$ be a finite group, then
every operator space $E$ satisfies the estimates
$\Constant_p^1(E,G), \,\ \Constant_{p'}^2(E,G) \le |G|^{1/p'}$ for
$1 \le p \le 2$.
\end{proposition}

\begin{proof} Let us assume that $\Constant_2^1(E,G) \le
|G|^{1/2}$ for every operator space $E$, then we have
$\Constant_2^2(E,G) = \Constant_2^1(E^{\star},G) \le |G|^{1/2}$ by
duality. The desired estimates are then obtained by complex
interpolation from the equalities $\Constant_1^1(E,G) =
\Constant_{\infty}^2(E,G)= 1$ --proved in \cite{GP}-- and the case
$p=2$. Therefore we focus our attention on the case $p = 2$. It
suffices to check that for all $m \ge 1$ and any family of
functions $\{f_{ij}: G \rightarrow E\}_{1 \le i,j \le m}$ $$\Big(
\sum_{\pi \in \G} d_{\pi} \Big\| \Big( \,\ \widehat{f}_{ij}(\pi)
\,\ \Big) \Big\|_{S_{d_{\pi}m}^2(E)}^2 \Big)^{1/2} \le |G|^{1/2}
\,\ \Big\| \Big( \,\ f_{ij} \,\ \Big) \Big\|_{S_m^2(L_E^2(G)).}$$
But if $G = \{g_1, g_2, \ldots g_n\}$, then
\begin{eqnarray*}
\Big\| \Big( \,\ \widehat{f}_{ij}(\pi) \,\ \Big)
\Big\|_{S_{d_{\pi}m}^2(E)} & = & \Big\| \Big( \,\ \frac{1}{n}
\sum_{k=1}^n f_{ij}(g_k) \pi(g_k)^{\star} \,\ \Big)
\Big\|_{S_{d_{\pi}m}^2(E)}
\\ & \le & \frac{1}{n} \sum_{k=1}^n \|\pi(g_k)^{\star}\|_{S_{d_{\pi}}^2}
\Big\| \Big( \,\ f_{ij}(g_k) \,\ \Big) \Big\|_{S_m^2(E)} \\ & \le
& d_{\pi}^{1/2} \,\ \Big\| \Big( \,\ f_{ij} \,\ \Big)
\Big\|_{S_m^2(L_E^2(G))}
\end{eqnarray*}
Therefore we obtain $$\Big\| \Big( \,\ \widehat{f}_{ij} \,\ \Big)
\Big\|_{S_m^2(\Lebesgue_E^2(\G))} \le \sqrt{\sum_{\pi \in \G}
d_{\pi}^2} \,\ \Big\| \Big( \,\ f_{ij} \,\ \Big)
\Big\|_{S_m^2(L_E^2(G))}$$ and, since $\displaystyle \sum_{\pi \in
\G} d_{\pi}^2 = |G|$ by the Peter-Weyl theorem, we are done.
\end{proof}

Next we show that we can not work with measure spaces $(\Omega,
\mathcal{A}, \nu)$ which are a \textbf{union of finitely many
$\nu$-atoms}. Before that we need to define the $cb$ distance
between two operator spaces. It is due to Pisier and it
constitutes the analog of the Banach-Mazur distance between two
Banach spaces in the context of operator space theory. Given two
operator spaces $E_1$ and $E_2$, we define their $cb$ distance by
the relation $d_{cb}(E_1,E_2) = \inf \{\|u\|_{cb(E_1,E_2)}
\|u^{-1}\|_{cb(E_2,E_1)}\}$ where the infimum runs over all
complete isomorphisms $u: E_1 \rightarrow E_2$. The following
result --also extracted from \cite{GP}-- relates the Fourier type
and cotype of two operator spaces $E_1$ and $E_2$ with their $cb$
distance.

\begin{lemma} \label{cb distance} Let $1 \le p \le 2$ and let
$E_1, E_2$ be operator spaces, then we have the estimates
$\Constant_p^1(E_2,G) \le d_{cb}(E_1,E_2) \,\
\Constant_p^1(E_1,G)$ and $\Constant_{p'}^2(E_2,G) \le
d_{cb}(E_1,E_2) \,\ \Constant_{p'}^2(E_1,G)$.
\end{lemma}

\begin{proposition} \label{Atoms} Let $1 \le p < q \le 2$ and
assume that $(\Omega, \mathcal{A}, \nu)$ is a union of finitely
many $\nu$-atoms. Then every compact group $G$ satisfies the
following estimates $$\begin{array}{rclcl}
\Constant_q^1(L^p(\Omega),G) & = &
\Constant_{q'}^2(L^{p'}(\Omega),G) & \le & \nu(\Omega)^{1/p - 1/q}
\\ \Constant_q^1(L^{p'}(\Omega),G) & = &
\Constant_{q'}^2(L^p(\Omega),G) & \le & \nu(\Omega)^{1/{q'} -
1/{p'}}. \end{array}$$ In fact, since $1/p - 1/q = 1/{q'} -
1/{p'}$, we have the same bound for both $cb$ norms.
\end{proposition}

\begin{proof} Applying lemma \ref{cb distance} and the last part
of theorem \ref{Lebesgue-Schatten} we get the estimates
$\Constant_q^1(L^p(\Omega), G) \le
d_{cb}(L^p(\Omega),L^q(\Omega))$ and
$\Constant_q^1(L^{p'}(\Omega),G) \le
d_{cb}(L^{p'}(\Omega),L^{q'}(\Omega))$. On the other hand it is
straightforward to check that, for $1 \le p_1 < p_2 \le \infty$
and such a measure space $(\Omega, \mathcal{A}, \nu)$, we have
$d_{cb}(L^{p_1}(\Omega),L^{p_2}(\Omega)) \le \nu(\Omega)^{1/p_1 -
1/p_2}$.
\end{proof}

In other words, we do not allow \textbf{finite dimensional}
Lebesgue spaces. Since the $cb$ distance between two Schatten
classes of the same finite dimension is also finite, the arguments
used in proposition \ref{Atoms} --theorem \ref{Lebesgue-Schatten}
and lemma \ref{cb distance}-- are also valid to show that the
unique Schatten classes with any possibility to make theorem
\ref{Lebesgue-Schatten} sharp are those of infinite dimension.

\subsection{Sufficient conditions}

The Fourier type and cotype of the subspaces of a given operator
space $E$ are bounded above by the respective type and cotype of
$E$. The proof of this result is straightforward, see \cite{GP}.

\begin{lemma} \label{Subspaces} Let $1 \le p \le 2$ and let $F$ be
a closed subspace of $E$, then we have the estimates
$\Constant_p^1(F,G) \le \Constant_p^1(E,G)$ and
$\Constant_{p'}^2(F,G) \le \Constant_{p'}^2(E,G)$.
\end{lemma}

After the conditions above, we shall work in the sequel with
infinite compact groups and infinite dimensional Lebesgue spaces
and Schatten classes. Since the measure space $(\Omega,
\mathcal{A}, \nu)$ is no longer a union of finitely many
$\nu$-atoms, we obtain that the $n$-dimensional space $l^p(n)$ is
a closed subspace of $L^p(\Omega)$ for all $n \ge 1$ and any $1
\le p \le \infty$. Moreover, recalling that the subspace of
diagonal matrices of $S_n^p$ is completely isomorphic to $l^p(n)$,
we deduce that the same happens for the Schatten classes
$S_{\N}^p$. Hence sharpness of theorem \ref{Lebesgue-Schatten}
will be guaranteed if, for $1 \le p < q \le 2$, we have
$$\begin{array}{llclcllc} (a') & \Constant_q^1(l^p(n),G) & = &
\Constant_{q'}^2(l^{p'}(n),G) & \longrightarrow & \infty & \quad
\mbox{as} \ \ n \rightarrow \infty & \qquad \qquad \\ (b') &
\Constant_q^1(l^{p'}(n),G) & = & \Constant_{q'}^2(l^p(n),G) &
\longrightarrow & \infty & \quad \mbox{as} \ \ n \rightarrow
\infty. &
\end{array}$$ Therefore our aim from now on will be the study of
the growth of the constants $\Constant_q^1(l^p(n),G)$ and
$\Constant_q^1(l^{p'}(n),G)$. The first remark about these
constants that we can already make is that both have a common
upper bound $$\Constant_q^1(l^p(n),G), \,\
\Constant_q^1(l^{p'}(n),G) \le n^{1/p - 1/q}.$$ This is an obvious
consequence of proposition \ref{Atoms}.

\subsection{Vector-valued spaces}

Theorem \ref{Lebesgue-Schatten} was also studied in \cite{GP} for
vector-valued spaces, here is the statement of the result
obtained.

\begin{theorem} \label{Vector-valued} Let $1 \le p \le \infty$ and
let $E$ be an operator space having Fourier type $\min(p,p')$
--respectively Fourier cotype $\max(p,p')$-- with respect to $G$.
Then the spaces $L_E^p(\Omega)$, $S_n^p(E)$ and $S_{\N}^p(E)$ have
Fourier type $\min(p,p')$ --respectively Fourier cotype
$\max(p,p')$-- with respect to $G$.
\end{theorem}

Let $1 \le p \le \infty$ and $\min(p,p') < q \le 2$. Let $E$ be as
in theorem \ref{Vector-valued}, then lemma \ref{Subspaces} gives
the following estimates $$\begin{array}{lclcrcr}
\Constant_q^1(L_E^p(\Omega),G) & \ge &
\Constant_q^1(L^p(\Omega),G)& \qquad \qquad
\Constant_q^1(L_E^p(\Omega),G) & \ge & \Constant_q^1(E,G) \\
\Constant_{q'}^2(L_E^p(\Omega),G) & \ge &
\Constant_{q'}^2(L^p(\Omega),G)& \qquad \qquad
\Constant_{q'}^2(L_E^p(\Omega),G) & \ge & \Constant_{q'}^2(E,G)
\end{array}$$ with the obvious modifications for the Schatten
classes. Hence we have shown that sharp Fourier type or cotype of
$L^p(\Omega)$ --respectively $S_{\N}^p$-- provides sharp Fourier
type or cotype of $L_E^p(\Omega)$ --respectively $S_{\N}^p(E)$--
with respect to $G$. Also the same conclusion is obtained assuming
sharp Fourier type or cotype of $E$. In particular, the sufficient
condition given above also works for vector-valued spaces.
Therefore we focus our attention on the growth of the constants
$\Constant_q^1(l^p(n),G)$ and $\Constant_q^1(l^{p'}(n),G)$.

\section{On the growth of $\Constant_q^1(l^p(n),G)$.}
\label{Section-p<2}

We shall assume in what follows that $G$ is a \textbf{compact
semisimple Lie group}. Semisimplicity is an essential assumption
in the arguments we shall be using. Anyway, for the moment, the
only property of such groups that we shall apply is the existence
of a maximal torus $\MT$ in $G$. The following result gives, in
particular, part $(a)$ in section \ref{Section-Statement} --with
the obvious modifications for Schatten classes-- whenever we work
with infinite dimensional operator spaces and compact semisimple
Lie groups.

\begin{theorem} \label{p<2}
Let $1 \le p < q \le 2$ and let $G$ be a compact semisimple Lie
group. Then there exists a constant $0 < \mathcal{K}(G,q) \le 1$
depending on $G$ and $q$ such that for all $n \ge 1$
$$\mathcal{K}(G,q) \,\ n^{1/p - 1/q} \le \Constant_q^1(l^p(n),G)
\le n^{1/p - 1/q}.$$
\end{theorem}

In particular we observe that the growth of
$\Constant_q^1(l^p(n),G)$ is optimal for compact semisimple Lie
groups. The proof of this result starts out applying the existence
of a maximal torus $\MT$ to consider a countable family $\{g_k: k
\ge 1\}$ of pairwise commuting elements of $G$, just take $g_k \in
\MT$. For every $n \ge 1$ we take $\Ball_n$ to be a neighborhood
of $\Identity$ --the identity element of $G$-- satisfying
$$g_j^{-1} \Ball_n \cap g_k^{-1} \Ball_n = \emptyset \qquad
\mbox{for} \ \ 1 \le j,k \le n \ \ \mbox{and} \ \ j \neq k.$$ We
recall here that we can always consider a central function $f_n$
supported in $\Ball_n$ and belonging to $L^q(G)$, for example take
$\Ball_n$ to be invariant under conjugations --see lemma $(5.24)$
of \cite{Fo}-- and $f_n = \Simple_{\Ball_n}$ where
$\Simple_{\Ball}$ stands for the cha\-racteristic function of
$\Ball$. Henceforth $f_n$ will be a central function in $L^q(G)$
supported in $\Ball_n$, to be fixed later. Then we define the
function $\Phi_n: G \rightarrow \C^n$ by $\Phi_n(g) = (f_n(g_1g),
f_n(g_2g), \ldots f_n(g_ng))$. We obviously have the estimate
$$\Constant_q^1(l^p(n),G) \ge
\frac{\|\widehat{\Phi}_n\|_{\Lebesgue_{l^p(n)}^{q'}(\G)}}
{\|\Phi_n\|_{L_{l^p(n)}^q(G)}}.$$ So it suffices to prove that
this quotient is bounded below by $\mathcal{K}(G,q) \,\ n^{1/p -
1/q}$. The following lemma will be very helpful for that purpose.

\begin{lemma} \label{Diagonal} Let $1 \le p_1, p_2 \le \infty$,
$\pi \in \G$ and $n \ge 1$. Consider the matrix-valued vector
$\mathrm{A}_{\pi,n} = (\pi(g_1), \pi(g_2), \ldots \pi(g_n))$, then
$$\|\mathrm{A}_{\pi,n}\|_{l_{S_{d_{\pi}}^{p_2}}^{p_1}(n)} =
\|\mathrm{A}_{\pi,n}\|_{S_{d_{\pi}}^{p_2}(l^{p_1}(n))} = n^{1/p_1}
d_{\pi}^{1/p_2}.$$
\end{lemma}

\begin{proof} Since $g_1, g_2, \ldots g_n$ are pairwise
commuting, there exists a basis of $\C^{d_{\pi}}$ of common
eigenvectors of $\pi(g_1), \pi(g_2) \ldots \pi(g_n)$. Therefore,
in that basis, all these matrices are diagonal $$\pi(g_k) = \left(
\begin{array}{ccc} \theta_1^k & & \\ & \ddots &
\\ & & \theta_{d_{\pi}}^k \end{array} \right).$$
Moreover, $|\theta_j^k| = 1$ for $1 \le j \le d_{\pi}$ because of
the unitarity of $\pi(g_k)$. Hence, applying the complete isometry
--see corollary $(1.3)$ of \cite{P2}-- between
$l_E^{p_2}(d_{\pi})$ and the subspace of diagonal matrices of
$S_{d_{\pi}}^{p_2}(E)$, we easily obtain the desired equality.
\end{proof}

\begin{itemize}

\item[$(i)$] \textbf{The value of
$\|\widehat{\Phi}_n\|_{\Lebesgue_{l^p(n)}^{q'}(\G)}$}. We begin by
recalling that, since $f_n$ is central, $$\widehat{f}_n (\pi) =
\frac{1}{d_{\pi}} \int_{G} f_n (g) \overline{\chi_{\pi}(g)} d
\mu(g) \,\ \mbox{{\Large 1}}_{d_{\pi}} = \gamma_{\pi,n}
\mbox{{\Large 1}}_{d_{\pi}}$$ by Schur's lemma. Here $\chi_{\pi}$
is the irreducible character associated to $\pi$ and
$\mbox{{\Large 1}}_m$ denotes the identity matrix of order $m
\times m$. On the other hand $f_n(g_k \cdot)$ is the translation
by $g_k$ of $f_n$, therefore $$\widehat{\Phi}_n(\pi) =
\frac{1}{d_{\pi}} \int_{G} f_n (g) \overline{\chi_{\pi}(g)} d
\mu(g) \,\ (\pi(g_1), \pi(g_2), \ldots \pi(g_n)) = \gamma_{\pi,n}
A_{\pi,n}.$$ So we get, by lemma \ref{Diagonal}, the following
equality
\begin{eqnarray*}
\|\widehat{\Phi}_n\|_{\Lebesgue_{l^p(n)}^{q'}(\G)} & = & \Big(
\sum_{\pi \in \G} d_{\pi} |\gamma_{\pi,n}|^{q'}
\|A_{\pi,n}\|_{S_{d_{\pi}}^{q'}(l^p(n))}^{q'} \Big)^{1/q'} \\ & =
& n^{1/p} \,\ \Big( \sum_{\pi \in \G} d_{\pi}^2
|\gamma_{\pi,n}|^{q'} \Big)^{1/q'} = n^{1/p} \,\
\|\widehat{f}_n\|_{\Lebesgue^{q'}(\G).}
\end{eqnarray*}

\item[$(ii)$] \textbf{The value of $\|\Phi_n\|_{L_{l^p(n)}^p(G)}$}.
We have
\begin{eqnarray*}
\|\Phi_n\|_{L_{l^p(n)}^q(G)} & = & \Big( \int_G \Big( \sum_{k=1}^n
|f_n(g_kg)|^p \Big)^{q/p} d \mu(g) \Big)^{1/q}
\\ & = & \Big( \sum_{k=1}^n \|f_n(g_k \cdot)\|_{L^q(G)}^q \Big)^{1/q} =
n^{1/q} \,\ \|f_n\|_{L^q(G)}
\end{eqnarray*}
since the sets $\{g_k^{-1} \Ball_n: 1 \le k \le n\}$ are pairwise
disjoint.
\end{itemize}

In summary, we have obtained that $\Constant_q^1(l^p(n),G) \ge
\mathcal{K}(G,q,n) \,\ n^{1/p - 1/q}$ where the constant
$\mathcal{K}(G,q,n)$ is given by $$\mathcal{K}(G,q,n) =
\frac{\|\widehat{f}_n\|_{\Lebesgue^{q'}(\G)}}
{\|f_n\|_{L^q(G)}}.$$ If we define $\mathcal{K}(G,q) = \inf_{n \ge
1} \mathcal{K}(G,q,n)$, it is obvious that $\mathcal{K}(G,q) \le
1$ by the Hausdorff-Young inequality on compact groups. Thus it
remains to check that $\mathcal{K}(G,q) > 0$. For that aim, since
we have not fixed $f_n$ yet, we need to see that $$\inf_{n \ge 1}
\sup \left\{
\frac{\|\widehat{f}\|_{\Lebesgue^{q'}(\G)}}{\|f\|_{L^q(G)}}: \,\ f
\ \ \mbox{central}, \ \ f \in L^q(G), \ \ \mbox{supp}(f) \subset
\Ball_n \right\}
> 0.$$ We shall prove this fact in section \ref{Section-Local}
where we study the supremum of the Hausdorff-Young quotient for
central functions supported in arbitrary small sets. As we shall
see immediately, semisimplicity of $G$ will be essential in our
proof.

\section{A simple expression for the Fourier transform \\ of central
functions} \label{Section-Simple}

In this section we apply some basic results concerning the
structure and representation theory of compact semisimple Lie
groups to provide a simple expression for the Fourier transform of
central functions defined on such groups. These algebraic
preliminaries can be found in Simon's book \cite{Si} or
alternatively in \cite{FH}, but we summarize here the main topics.
Let $G$ be a compact semisimple Lie group and let $\g$ be its Lie
algebra. In what follows we choose once and for all an explicit
maximal torus $\MT$ in $G$ while $\h$ will stand for its Lie
algebra. That is, $\h$ is the Cartan subalgebra of $\g$. The rank
of $G$ will be denoted by $r$, in particular $\MT \simeq \T^r$
where $\T = \R / \Z$ with its natural group structure. Also, as it
is customary, we consider the complexification $\g_{\C} = \g
\oplus i \g$ --with complex conjugates taken so that $\g_{\R} =
\{Z \in \g_{\C}: Z = \overline{Z}\} = i \g$ and similarly $\h_{\R}
= i \h$-- with the complex inner product $\langle \,\ , \,\
\rangle$ induced by the Killing form. We also recall that the Weyl
group $\Weyl$ associated to $G$ can be seen as a set of $r \times
r$ unitary matrices $W$ --isometries on $\h_{\R}$-- with integer
entries and $\det W = \pm 1$. In particular the set $\Weyl^{\star}
= \{W^t: W \in \Weyl\}$ becomes a set of isometries on
$\h_{\R}^{\star}$. The symbol $\Root$ will stand for the set of
roots while, if we take $H_0 \in \h_{\R}$ such that $\alpha(H_0)
\neq 0$ for any root $\alpha$, the symbol $\Root^+ = \{\alpha \in
\Root: \alpha(H_0) > 0\}$ denotes the set of positive roots.
Finally we shall write $\Lambda_{\mathrm{W}}$ and
$\Lambda_{\mathrm{DW}}$ for the weight lattice and the set of
dominant weights respectively.

Once we have fixed some notation, let us consider a central
function $f: G \rightarrow \C$ and a dominant weight $\lambda \in
\Lambda_{\mathrm{DW}}$. By the dominant weight theorem there
exists a unique $\pi_{\lambda} \in \G$ associated to $\lambda$
and, since $f$ is central, we can write by Schur's lemma
$$\widehat{f}(\pi_{\lambda}) = \frac{1}{d_{\lambda}} \int_G f(g)
\overline{\chi_{\lambda}(g)} d \mu(g) \,\ \mbox{{\Large
1}}_{d_{\lambda}}$$ where $d_{\lambda}$ is the degree of
$\pi_{\lambda}$, $\chi_{\lambda}$ is the character of
$\pi_{\lambda}$ and $\mbox{{\Large 1}}_m$ denotes the $m \times m$
identity matrix. We now recall the definition of the functions
$A_{\beta}$ appearing in the Weyl character formula. Given $\beta
\in \h_{\R}^{\star}$, we define the functions $\exp_{\beta}:
\h_{\R} \rightarrow \C$ and $A_{\beta}: \h_{\R} \rightarrow \C$ by
the relations
\begin{eqnarray*}
\exp_{\beta}(H) & = & e^{2 \pi i \langle \beta,H \rangle} \\
A_{\beta}(H) & = & \sum_{W \in \Weyl} \det W \,\ \exp_{\beta}\big(
W(H) \big).
\end{eqnarray*}
The maximal torus $\MT$ is isomorphic via the exponential mapping
to the quotient space $\h_{\R} / \mathrm{L}_{\mathrm{W}}$, where
$\mathrm{L}_{\mathrm{W}}$ is the set of those $H \in \h_{\R}$
satisfying $\exp(2 \pi i H) = \Identity$. That is,
$\mathrm{L}_{\mathrm{W}}$ is the dual lattice of
$\Lambda_{\mathrm{W}}$. Therefore, the functions $\exp_{\beta}$
and $A_{\beta}$ are well-defined functions on $\MT$ if and only if
$\beta \in \Lambda_{\mathrm{W}}$. As it is well known, the
integral form $$\delta = \frac{1}{2} \sum_{\alpha \in \Root^{+}}
\alpha$$ is not necessarily a weight and so the functions
$\exp_{\delta}$ and $A_{\delta}$ could be not well-defined on
$\MT$. To avoid this difficulty we assume for the moment that $G$
is \textbf{simply connected}, this condition on $G$ assures that
$\delta \in \Lambda_{\mathrm{W}}$. Hence, applying consecutively
the Weyl integration formula and the Weyl character formula, we
get
\begin{eqnarray*}
\widehat{f}(\pi_{\lambda}) & = & \frac{1}{d_{\lambda} |\Weyl|}
\int_{\MT} f(t) \overline{\chi_{\lambda}(t)} \,\ |A_{\delta}(t)|^2
d m(t) \,\ \mbox{{\Large 1}}_{d_{\lambda}} \\ & = &
\frac{1}{d_{\lambda} |\Weyl|} \int_{\MT} f(t) A_{\delta}(t) \,\
\overline{A_{\lambda + \delta}(t)} d m(t) \,\ \mbox{{\Large
1}}_{d_{\lambda}}
\end{eqnarray*}
where $m$ denotes the Haar measure on $\MT$ normalized so that
$m(\MT) = 1$. Now, if we write $A_{\lambda + \delta}$ as a linear
combination of exponentials, we obtain
\begin{eqnarray*}
\widehat{f}(\pi_{\lambda}) & = & \frac{1}{d_{\lambda} |\Weyl|}
\sum_{W \in \Weyl} \det W \int_{\MT} f(t) A_{\delta}(t) \,\
\exp_{-(\lambda + \delta)} (W(t)) d m(t) \,\ \mbox{{\Large
1}}_{d_{\lambda}} \\ & = & \frac{1}{d_{\lambda}} \int_{\MT} f(t)
A_{\delta}(t) \,\ \exp_{-(\lambda + \delta)}(t) d m(t) \,\
\mbox{{\Large 1}}_{d_{\lambda}}
\end{eqnarray*}
since $A_{\delta}(W(t)) = \det W A_{\delta}(t)$ and $f(W(t)) =
f(t)$. We recall that, taking coordinates with respect to the
basis $\{\omega_1, \omega_2, \ldots \omega_r\}$ of fundamental
weights, any weight $\lambda \in \Lambda_{\mathrm{W}}$ has integer
coordinates. Therefore, we can understand the last expression as
the Fourier transform of $f A_{\delta}$ on the maximal torus $\MT$
evaluated at $\lambda + \delta$. Hence we have
\begin{equation} \label{1}
\widehat{f}(\pi_{\lambda}) = \frac{1}{d_{\lambda}} \,\
\mathcal{F}_{\MT} (f A_{\delta}) (\lambda + \delta) \ \
\mbox{{\Large 1}}_{d_{\lambda}}
\end{equation}
for $f: G \rightarrow \C$ central and $G$ any compact semisimple
simply connected Lie group. When $G$ is not simply connected, a
more careful approach is needed. We have $W^t(\delta) \pm \delta
\in \Lambda_{\mathrm{W}}$ for all $W \in \Weyl$. In particular we
note that $$\exp_{\pm \delta} A_{\lambda + \delta} = \sum_{W \in
\Weyl} \det W \exp_{W^t(\lambda + \delta) \pm \delta}$$ is a
well-defined function on $\MT$ for all $\lambda \in
\Lambda_{\mathrm{DW}}$. This remark allows us to write
$\overline{\chi_{\lambda}} \,\ |A_{\delta}|^2 = (\exp_{\delta}
\overline{A_{\lambda + \delta}}) \,\ (\exp_{- \delta} A_{\delta})$
as a well-defined function on $\MT$. Henceforth, applying again
Schur's lemma, the Weyl integration formula and the Weyl character
formula, we get
\begin{eqnarray*}
\widehat{f}(\pi_{\lambda}) & = & \frac{1}{d_{\lambda} |\Weyl|}
\sum_{W \in \Weyl} \det W \int_{\MT} f(t) (\exp_{- \delta}
A_{\delta}) (t) \,\ \exp_{\delta - W^t(\lambda + \delta)} (t) d
m(t) \,\ \mbox{{\Large 1}}_{d_{\lambda}} \\ & = &
\frac{1}{d_{\lambda}} \int_{\MT} f(t) (\exp_{- \delta} A_{\delta})
(t) \,\ \exp_{- \lambda} (t) dm(t) \,\ \mbox{{\Large
1}}_{d_{\lambda}}
\end{eqnarray*}
where the last equality follows from the change of variable $t
\mapsto W^t(t)$. That is, we have shown that
\begin{equation} \label{2}
\widehat{f}(\pi_{\lambda}) = \frac{1}{d_{\lambda}} \,\
\Fourier_{\MT}(f B_{\delta})(\lambda) \,\ \mbox{{\Large
1}}_{d_{\lambda}}
\end{equation}
where $B_{\delta} = \exp_{- \delta} A_{\delta}$. This expression
is now valid for any compact semisimple Lie group and it coincides
with $(\ref{1})$ for simply connected ones.

\section{A local variant of the Hausdorff-Young inequality \\ on
compact semisimple Lie groups} \label{Section-Local}

As we mentioned in the introduction, this section is devoted to
the proof of a local variant of the Hausdorff-Young inequality on
compact semisimple Lie groups. We recall that this result provides
the relation $\mathcal{K}(G,q) > 0$ for $1 \le q \le 2$, that we
needed in section \ref{Section-p<2}.

\begin{theorem} \label{Local}
Let $1 \le q \le 2$ and let $G$ be a compact semisimple Lie group.
Then there exists a constant $0 < \mathcal{K}(G,q) \le 1$ such
that, for any open set $\Ball \subset G$, we have $$\sup \left\{
\frac{\|\widehat{f}\|_{\Lebesgue^{q'}(\G)}}{\|f\|_{L^q(G)}}: \ \ f
\ \ \mbox{central}, \ \ f \in L^q(G), \ \ \textnormal{supp}(f)
\subset \Ball \right\} \ge \mathcal{K}(G,q).$$
\end{theorem}

Since the norms of $\widehat{f}$ and $f$ --on $\Lebesgue^{q'}(\G)$
and $L^q(G)$ respectively-- do not change under translations of
$f$, we can assume without loss of generality that $\Ball$ is a
neighborhood of $\Identity$. Before the proof of theorem
\ref{Local} we need some auxiliary results. Let us assume that $G$
is simply connected and let $f: G \rightarrow \C$ be a central
function. A quick look at relation $(\ref{1})$ given above, allows
us to write
\begin{equation} \label{3}
\widehat{f}(\pi_{\lambda}) = \frac{1}{d_{\lambda}} \,\ \det W \,\
\mathcal{F}_{\MT} (f A_{\delta}) (W^t(\lambda + \delta)) \ \
\mbox{{\Large 1}}_{d_{\lambda}}
\end{equation}
for all $W \in \Weyl$. On the other hand, let us denote by
$P_{\alpha}$ the hyperplane of $\h_{\R}^{\star}$ orthogonal to
$\alpha$ with respect to the complex inner product given by the
Killing form. The infinitesimal Cartan-Stiefel diagram is then
given by the expression $$\Planes = \bigcup_{\alpha \in \Root}
P_{\alpha}.$$

\begin{lemma} \label{Bijection}
Let $G$ be a compact semisimple simply connected Lie group. Then
we have $\{W^t(\lambda + \delta): \,\ W \in \Weyl, \,\ \lambda \in
\Lambda_{\mathrm{DW}}\} = \Lambda_{\mathrm{W}} \setminus \Planes$.
Moreover, the mapping $(W, \lambda) \in \Weyl \times
\Lambda_{\mathrm{DW}} \mapsto W^t(\lambda + \delta) \in
\Lambda_{\mathrm{W}} \setminus \Planes$ is injective.
\end{lemma}

\begin{proof} Since $G$ is simply connected we have that $\{
\lambda + \delta: \lambda \in \Lambda_{\mathrm{DW}}\} =
\Lambda_{\mathrm{W}} \cap \displaystyle
\Chamber^{\mbox{{\footnotesize int}}}$. Here $\Chamber$ stands for
the fundamental Weyl chamber and $\Chamber^{\mbox{{\footnotesize
int}}}$ for its interior. Now, since $\Planes$ and
$\Lambda_{\mathrm{W}}$ are invariant under the action of
$\Weyl^{\star}$ and for any Weyl chamber $\mathrm{C}$ there exists
a unique $W \in \Weyl$ such that $W^t(\Chamber) = \mathrm{C}$, we
obtain the desired equality. Finally, the injectivity follows from
the uniqueness mentioned above.
\end{proof}

\begin{proposition} \label{Norm}
Let $G$ be a compact semisimple simply connected Lie group and let
$f: G \rightarrow \C$ be a central function. Then there exists a
constant $\mathcal{A}(G,q)$ depending on $G$ and $q$, such that
$$\|\widehat{f}\|_{\Lebesgue^{q'}(\G)} = \mathcal{A}(G,q) \,\
\Big[ \sum_{\lambda \in \Lambda_{\mathrm{W}} \setminus \Planes}
\frac{|\Fourier_{\MT}(f A_{\delta})(\lambda)|^{q'}}{\displaystyle
\prod_{\alpha \in \Root^{+}} |\langle \alpha, \lambda \rangle|^{q'
- 2}} \Big]^{1/q'}.$$
\end{proposition}

\begin{proof} Since $f$ is central and $G$ is simply connected, we
can apply expression $(\ref{3})$ to obtain
\begin{eqnarray*}
\|\widehat{f}\|_{\Lebesgue^{q'}(\G)} & = & \Big[ \sum_{\lambda \in
\Lambda_{\mathrm{DW}}} d_{\lambda} \|\widehat{f}
(\pi_{\lambda})\|_{S^{q'}_{d_{\lambda}}}^{q'} \Big]^{1/q'} \\ & =
& \Big[ \frac{1}{|\Weyl|} \sum_{W \in \Weyl} \sum_{\lambda \in
\Lambda_{\mathrm{DW}}} d_{\lambda} \Big| \frac{1}{d_{\lambda}}
\Fourier_{\MT} (f A_{\delta}) (W^t(\lambda + \delta)) \Big|^{q'}
\|\mbox{{\Large 1}}_{d_{\lambda}}\|_{S^{q'}_{d_{\lambda}}}^{q'}
\Big]^{1/q'}
\end{eqnarray*}
Moreover, the Weyl dimension formula for $d_{\lambda}$ gives
$$\|\widehat{f}\|_{\Lebesgue^{q'}(\G)} = \Big[ \frac{1}{|\Weyl|}
\prod_{\alpha \in \Root^{+}} |\langle \alpha, \delta
\rangle|^{q'-2} \sum_{W \in \Weyl} \sum_{\lambda \in
\Lambda_{\mathrm{DW}}} \frac{|\Fourier_{\MT} (f A_{\delta})
(W^t(\lambda + \delta))|^{q'}}{\displaystyle \prod_{\alpha \in
\Root^{+}} |\langle \alpha, \lambda + \delta \rangle|^{q'-2}}
\Big]^{1/q'}.$$ Finally we observe that $$\prod_{\alpha \in
\Root^{+}} |\langle \alpha, \lambda + \delta \rangle| =
\prod_{\alpha \in \Root} |\langle W(\alpha), \lambda + \delta
\rangle|^{1/2} = \prod_{\alpha \in \Root^{+}} |\langle \alpha,
W^t(\lambda + \delta) \rangle|$$ since any $W \in \Weyl$ is a
permutation of the set of roots. Therefore, by lemma
\ref{Bijection} we have $$\|\widehat{f}\|_{\Lebesgue^{q'}(\G)} =
\Big[ \frac{1}{|\Weyl|} \prod_{\alpha \in \Root^{+}} |\langle
\alpha, \delta \rangle|^{q'-2} \sum_{\lambda \in
\Lambda_{\mathrm{W}} \setminus \Planes} \frac{|\Fourier_{\MT} (f
A_{\delta}) (\lambda)|^{q'}}{\displaystyle \prod_{\alpha \in
\Root^{+}} |\langle \alpha, \lambda \rangle|^{q'-2}}
\Big]^{1/q'}.$$ The proof is completed just by taking
$\displaystyle \mathcal{A}(G,q) = \Big[ \frac{1}{|\Weyl|}
\prod_{\alpha \in \Root^{+}} |\langle \alpha, \delta
\rangle|^{q'-2} \Big]^{1/q'}$.
\end{proof}

We are now ready to give the proof of theorem \ref{Local} for
simply connected groups. Let $\{H_1, H_2, \ldots H_r\}$ be the
predual basis of the fundamental weights, any element of
$\mathrm{L}_{\mathrm{W}}$ can be written as a linear combination
of $H_1, H_2, \ldots H_r$ with integer coefficients. Then, since
$\MT \simeq \h_{\R} / \mathrm{L}_{\mathrm{W}}$, we can regard
$\MT$ as the subset of $\h_{\R}$ $$\mt = \Big\{ \sum_{k=1}^r x_k
H_k: \,\ - 1/2 \le x_k < 1/2 \Big\}.$$ On the other hand, let us
fix a bounded central function $f_0: G \rightarrow \C$, then $f_0$
can be understood as a function on $\MT$ invariant under the
action of $\Weyl$. Now, since the Weyl group is generated by a set
of reflections in $\h_{\R}$, $f_0$ can be regarded as a
complex-valued function on $\h_{\R}$, supported in $\mt$ and
symmetric under such reflections. Let us recall that $\{\omega_1,
\omega_2, \ldots \omega_r\}$ stands for the basis of fundamental
weights. Let $\tau = 1 - 2/q'$, the way we have interpreted the
function $f_0$ allows us to define the function
$$\widehat{I_{\tau}(f_0 A_{\delta})}: \h_{\R}^{\star}
\longrightarrow \C \qquad \mbox{as}$$ $$\widehat{I_{\tau}(f_0
A_{\delta})} (\xi) = \frac{1}{\displaystyle \prod_{\alpha \in
\Root^{+}} |\langle \alpha, \xi \rangle|^{\tau}} \,\
\Fourier_{\h_{\R}} (f_0 A_{\delta}) (\xi) \qquad \mbox{where}
\quad \xi = \sum_{k=1}^r \xi_k \omega_k.$$

\begin{remark}
The motivation for the notation employed is that in a classical
group such as $SU(2)$ the function just defined is nothing but the
Fourier transform of the fractional integral operator
$$I_{\tau}(f) (x) = \frac{1}{\Gamma(\tau)} \int_{- \infty}^x f(y)
(x-y)^{\tau-1} dy$$ acting on $f_0 A_{\delta}$. Here lies the main
difference with the commutative case --where a Hausdorff-Young
inequality of local type has been already investigated, see
\cite{An}-- since the presence of the degrees $d_{\lambda}$ --as a
product in proposition \ref{Norm} by the Weyl dimension formula--
requires the presence of a factor of $\Fourier_{\h_{\R}}(f_0
A_{\delta})$. This fact does not happen in the commutative case
since $d_{\lambda} = 1$ for all $\lambda \in
\Lambda_{\mathrm{DW}}$.
\end{remark}

\begin{lemma} \label{Zeros}
Let $G$ be a compact semisimple simply connected Lie group and let
$f: G \rightarrow \C$ be a central function. Then we have
$\Fourier_{\h_{\R}} (f A_{\delta}) (\xi) = 0$ for all $\xi \in
\Planes$.
\end{lemma}

\begin{proof}
If $\xi \in \Planes$, there exists a root $\alpha$ such that $\xi
\in P_{\alpha}$. Let $S_{\alpha}$ be the reflection in
$P_{\alpha}$, then $\Fourier_{\h_{\R}} (f A_{\delta}) (\xi) = \det
S_{\alpha} \,\ \Fourier_{\h_{\R}} (f A_{\delta}) (S_{\alpha}(\xi))
= - \Fourier_{\h_{\R}} (f A_{\delta}) (\xi)$ since, as it is well
known, $S_{\alpha} \in \Weyl^{\star}$.
\end{proof}

The function $\Fourier_{\h_{\R}}(f_0 A_{\delta})$ is analytic
since $f_0 A_{\delta}$ has compact support and, by lemma
\ref{Zeros}, it vanishes at $$\Planes = \{\xi \in \h_{\R}^{\star}:
\prod_{\alpha \in \Root^+} \langle \alpha, \xi \rangle = 0\}.$$ In
particular, since $0 \le \tau < 1$, $\widehat{I_{\tau}(f_0
A_{\delta})}$ is continuous and takes the value $0$ on $\Planes$.
Now we write the norm of this function in terms of a Riemann sum
$$\|\widehat{I_{\tau}(f_0 A_{\delta})}\|_{L^{q'}(\h_{\R}^{\star})}
= \lim_{k \rightarrow \infty} \Big[ \sum_{\lambda \in
\Lambda_{\mathrm{W}}} \frac{\mathrm{V}_G}{k^r} \,\
\frac{|\Fourier_{\h_{\R}} (f_0 A_{\delta}) (k^{-1}
\lambda)|^{q'}}{\displaystyle \prod_{\alpha \in \Root^{+}}
|\langle \alpha, k^{-1} \lambda \rangle|^{\tau q'}}
\Big]^{1/q'},$$ where $\mathrm{V}_G$ denotes the volume of a cell
of $\Lambda_{\mathrm{W}}$. Moreover $\phi_k(x) = k^{\sigma}
f_0(kx) A_{\delta}(kx)$ is supported in $\mt$ and the relation
$\Fourier_{\h_{\R}} (f_0 A_{\delta}) (k^{-1} \lambda) = k^{r -
\sigma} \Fourier_{\MT} (\phi_k) (\lambda)$ is satisfied for all
$\lambda \in \Lambda_{\mathrm{W}}$. Taking $\sigma = \tau
|\Root^{+}| + r/q$, we obtain $$\|\widehat{I_{\tau}(f_0
A_{\delta})}\|_{L^{q'}(\h_{\R}^{\star})} = \mathrm{V}_G^{1/q'}
\lim_{k \rightarrow \infty} \Big[ \sum_{\lambda \in
\Lambda_{\mathrm{W}} \setminus \Planes} \frac{|\Fourier_{\MT}
(\phi_k) (\lambda)|^{q'}}{\displaystyle \prod_{\alpha \in
\Root^{+}} |\langle \alpha, \lambda \rangle|^{\tau q'}}
\Big]^{1/q'},$$ since we know that for $\lambda \in \Planes$ we
get nothing. Finally, let us define $\varphi_k: \h_{\R}
\rightarrow \C$ by the relation $\phi_k = \varphi_k A_{\delta}$.
The function $\varphi_k$ satisfies $\varphi_k(W(x))= \varphi_k(x)$
for all $W \in \Weyl$ and is supported in $k^{-1} \mt$, hence we
can understand $\varphi_k$ as a central function on $G$. We can
also say that, as a consequence of the well known relation
\begin{equation} \label{4}
A_{\delta} = \exp_{- \delta} \prod_{\alpha \in \Root^{+}}
(\exp_{\alpha} - 1),
\end{equation}
$\varphi_k$ has no singularities. Therefore proposition \ref{Norm}
provides the following relation for some constant
$\mathcal{B}(G,q)$ depending on $G$ and $q$
\begin{equation} \label{5}
\|\widehat{I_{\tau}(f_0 A_{\delta})}\|_{L^{q'}(\h_{\R}^{\star})} =
\mathcal{B}(G,q) \,\ \lim_{k \rightarrow \infty}
\|\widehat{\varphi}_k\|_{\Lebesgue^{q'}(\G)}.
\end{equation}

On the other hand, since $\varphi_k$ can be seen as a central
function on $G$, we can estimate the norm of $\varphi_k$ on
$L^q(G)$. By the Weyl integration formula we get
\begin{eqnarray*}
\|\varphi_k\|_{L^q(G)} & = & \Big[ \frac{1}{|\Weyl|} \int_{\MT}
|\varphi_k A_{\delta} (t)|^q |A_{\delta}(t)|^{2-q} d m(t)
\Big]^{1/q} \\ & = & \Big[ \frac{k^{\sigma q}}{|\Weyl|} \int_{\mt}
|f_0 A_{\delta} (kx)|^q |A_{\delta}(x)|^{2-q} dx \Big]^{1/q} \\ &
\le & \Big[ \frac{(2 \pi)^{(2-q) |\Root^{+}|}}{|\Weyl|} \,\
k^{\sigma q} \int_{\mt} |f_0 A_{\delta} (kx)|^q \prod_{\alpha \in
\Root^{+}} |\alpha(x)|^{2-q} dx \Big]^{1/q},
\end{eqnarray*}
where the last inequality follows from $(\ref{4})$. Now, under the
change of variable $y = kx$ and taking $\mathcal{C}(G,q) = (2
\pi)^{\tau |\Root^{+}|} |\Weyl|^{-1/q}$, we obtain
$$\|\varphi_k\|_{L^q(G)} \le \mathcal{C}(G,q) \,\ k^{\sigma - \tau
|\Root^{+}| - r/q} \,\ \Big( \int_{\mt} |f_0 A_{\delta} (y)|^q
\prod_{\alpha \in \Root^{+}} |\alpha(y)|^{\tau q} dy
\Big)^{1/q}.$$ Recall that $\mbox{supp} (f_0 A_{\delta}) \subset
\mt$, therefore the integral over $k \mt$ --the domain of
integration after the change of variable-- reduces to the same
integral over $\mt$. But $\sigma - \tau |\Root^{+}| - r/q = 0$ and
the product inside the integral is bounded over $\mt$, say by
$\mathrm{M}_G$. Then we can write
\begin{equation} \label{6}
\|\varphi_k\|_{L^q(G)} \le \mathcal{C}(G,q) \,\ \mathrm{M}_G \,\
\|f_0 A_{\delta}\|_{L^q(\h_{\R})}.
\end{equation}

In summary, by $(\ref{5})$ and $(\ref{6})$, we know there exists a
constant $\mathcal{D}(G,q)$ depending on $G$ and $q$ such that
$$\mathcal{K}(G,q) = \mathcal{D}(G,q) \,\
\frac{\|\widehat{I_{\tau}(f_0
A_{\delta})}\|_{L^{q'}(\h_{\R}^{\star})}}{\|f_0
A_{\delta}\|_{L^q(\h_{\R})}} \le \liminf_{k \rightarrow \infty}
\frac{\|\widehat{\varphi}_k\|_{\Lebesgue^{q'}(\G)}}{\|\varphi_k\|_{L^q(G)}}
\le 1.$$ Since $f_0$ is bounded we easily obtain that $f_0
A_{\delta} \in L^q(\h_{\R})$, $\widehat{I_{\tau}(f_0 A_{\delta})}
\in L^{q'}(\h_{\R}^{\star})$ and $\mathcal{K}(G,q) > 0$. Therefore
we have found a family $\{\varphi_k: \,\ k \ge 1\}$ of central
functions on $G$ whose supports are eventually in $\Ball$ and such
that their Hausdorff-Young quotient of exponent $q$ is bounded
below by a positive constant. This concludes the proof of theorem
\ref{Local} for compact semisimple simply connected Lie groups.

If $G$ is not simply connected, some extra comments have to be
made. In any case we shall not give complete proofs of any of
them, the details are left to the reader.

\begin{itemize}
\item[$(i)$] Generalization $(\ref{3}$) of formula $(\ref{1})$ has
no meaning here, but we can generalize formula $(\ref{2})$ as
$$\widehat{f}(\pi_{\lambda}) = \frac{1}{d_{\lambda}} \,\ \det W
\,\ \Fourier_{\MT} (f B_{\delta}) (W^t(\lambda + \delta) - \delta)
\,\ \mbox{{\Large 1}}_{d_{\lambda}}.$$ This generalization
provides a couple of results parallel to lemmas \ref{Bijection}
and \ref{Zeros}. Namely,
\begin{itemize}
\item We have $\{W^t(\lambda + \delta) - \delta: \,\ W \in \Weyl, \,\
\lambda \in \Lambda_{\mathrm{DW}}\} = \Lambda_{\mathrm{W}}
\setminus (\Planes - \delta)$. The mapping $(W, \lambda) \in \Weyl
\times \Lambda_{\mathrm{DW}} \mapsto W^t(\lambda + \delta) -
\delta \in \Lambda_{\mathrm{W}} \setminus (\Planes - \delta)$ is
injective.
\item If $f: G \rightarrow \C$ is central, then
$\Fourier_{\h_{\R}}(f B_{\delta})(\xi) = 0$ for all $\xi \in
\Planes - \delta$.
\end{itemize}
\item[$(ii)$] Proposition \ref{Norm} is now replaced by the
following identity, valid for central functions $f: G \rightarrow
\C$ $$\|\widehat{f}\|_{\Lebesgue^{q'}(\G)} = \mathcal{A}(G,q) \,\
\Big[ \sum_{\lambda \in \Lambda_{\mathrm{W}} \setminus (\Planes -
\delta)} \frac{|\Fourier_{\MT}(f
B_{\delta})(\lambda)|^{q'}}{\displaystyle \prod_{\alpha \in
\Root^{+}} |\langle \alpha, \lambda + \delta \rangle|^{q' - 2}}
\Big]^{1/q'}.$$
\item[$(iii)$] The bases of $\h_{\R}^{\star}$ and $\h_{\R}$
respectively which generate $\Lambda_{\mathrm{W}}$ and
$\mathrm{L}_{\mathrm{W}}$ with integer coefficients are no longer
the basis of fundamental weights and its predual. In fact, the
fundamental weights generate the weight lattice of the universal
covering group of $G$, which is a lattice containing
$\Lambda_{\mathrm{W}}$ and strictly bigger than it. Therefore we
need to define $\{H_1, H_2, \ldots H_r\}$ and $\{\omega_1,
\omega_2, \ldots \omega_r\}$ just as the bases --of $\h_{\R}$ and
$\h_{\R}^{\star}$ respectively-- for which
$\mathrm{L}_{\mathrm{W}}$ and $\Lambda_{\mathrm{W}}$ have integer
coefficients. Once we have clarified this point, we can define
$\mt$ in the same way and regard $f_0$ as a bounded complex-valued
function on $\h_{\R}$, supported in $\mt$ and symmetric under the
reflections that generate $\Weyl$.
\item[$(iv)$] Let us recall that if $\delta \notin
\Lambda_{\mathcal{W}}$, the function $A_{\delta}$ is not
well-defined on $\MT$. But $A_{\delta}$ is originally defined on
$\h_{\R}$ and $\delta \notin \Lambda_{\mathcal{W}}$ is not an
obstacle to work with $A_{\delta}$ as a function defined on
$\h_{\R}$. On the other hand, $(ii)$ leads us to consider --in the
same spirit as in the proof given for simply connected groups--
the function $$\widehat{\widetilde{I}_{\tau} (f_0 B_{\delta})}
(\xi) = \frac{1}{\displaystyle \prod_{\alpha \in \Root^{+}}
|\langle \alpha, \xi + \delta \rangle|^{\tau}} \,\
\Fourier_{\h_{\R}} (f_0 B_{\delta}) (\xi).$$ Now, the remark given
about $A_{\delta}$ shows that $\widehat{\widetilde{I}_{\tau} (f_0
B_{\delta})} (\xi) = \widehat{I_{\tau} (f_0 A_{\delta})} (\xi +
\delta)$. Hence we can proceed as before expressing the norm of
this function in $L^{q'}(\h_{\R}^{\star})$ as a Riemann sum, but
this time we take the lattice $\Lambda_{\mathrm{W}} + \delta$
instead of $\Lambda_{\mathrm{W}}$
$$\|\widehat{\widetilde{I}_{\tau} (f_0
B_{\delta})}\|_{L^{q'}(\h_{\R}^{\star})} = \lim_{k \rightarrow
\infty} \Big[ \sum_{\lambda \in \Lambda_{\mathrm{W}} + \delta}
\frac{\mathrm{V}_G}{k^r} \,\ \frac{|\Fourier_{\h_{\R}} (f_0
A_{\delta}) (k^{-1} \lambda)|^{q'}}{\displaystyle \prod_{\alpha
\in \Root^{+}} |\langle \alpha, k^{-1} \lambda \rangle|^{\tau q'}}
\Big]^{1/q'}.$$
\item[$(v)$] It is not difficult to check that
$\Fourier_{\h_{\R}} (f_0 A_{\delta}) (k^{-1} \lambda) = k^{r -
\sigma} \Fourier_{\MT} (\varphi_k B_{\delta}) (\lambda - \delta)$,
where $\varphi_k$ is defined as we did above. Hence we get
\begin{eqnarray*}
\|\widehat{\widetilde{I}_{\tau} (f_0
B_{\delta})}\|_{L^{q'}(\h_{\R}^{\star})} & = & \mathrm{V}_G^{1/q'}
\lim_{k \rightarrow \infty} \Big[ \sum_{\lambda \in
\Lambda_{\mathrm{W}} \setminus (\Planes - \delta)}
\frac{|\Fourier_{\MT} (\varphi_k B_{\delta})
(\lambda)|^{q'}}{\displaystyle \prod_{\alpha \in \Root^{+}}
|\langle \alpha, \lambda + \delta \rangle|^{\tau q'}} \Big]^{1/q'}
\\ & = & \mathcal{B}(G,q) \,\ \lim_{k \rightarrow \infty}
\|\widehat{\varphi}_k\|_{\Lebesgue^{q'}(\G)}.
\end{eqnarray*}
\end{itemize}
Finally, to estimate the norm of $\varphi_k$ on $L^q(G)$, we
follow the same arguments. This completes the proof of theorem
\ref{Local} and, consequently, the proof of theorem \ref{p<2}.

\begin{remark}
Let $\{\Ball_n: n \ge 1\}$ be a basis of neighborhoods of
$\Identity$, and let $$\mathcal{K}(G,q) = \inf_{n\ge 1} \sup
\left\{
\frac{\|\widehat{f}\|_{\Lebesgue^{q'}(\G)}}{\|f\|_{L^q(G)}}: f \ \
\mbox{central}, \ \ f \in L^q(G), \ \ \mbox{supp}(f) \subset
\Ball_n \right\}.$$ This constant does not depend on the chosen
basis and theorem \ref{Local} states that $ 0 < \mathcal{K}(G,q)
\le 1$ for any $1 \le q \le 2$ and any compact semisimple Lie
group. However, it would be interesting to find the exact value of
that constant. Sharp constants for the Hausdorff-Young inequality
were investigated in \cite{Ba}, \cite{B1} or \cite{Ru}. In the
local case, if $\mathcal{B}_q = \sqrt{q^{1/q}/q^{'1/q'}}$ stands
for the Babenko-Beckner constant, it is already known that
$\mathcal{K}(\T,q) = \mathcal{B}_q$. Andersson proved it for $q'$
an even integer in \cite{An} and Sj\"{o}lin completed the proof, see
\cite{Sj}. Also it is obvious that $\mathcal{K}(G,1) =
\mathcal{K}(G,2) = 1$ for any compact group $G$. In the general
case, a detailed look at the proof of theorem \ref{Local} gives
that the constant $\mathcal{K}(G,q)$ is the supremum of
$$|\Weyl|^{\tau} \prod_{\alpha \in \Root^+} |\langle \alpha,
\delta \rangle|^{\tau} \,\ \mathrm{V}_G^{-1/q'} \,\ \lim_{k
\rightarrow \infty} \frac{\displaystyle \Big(
\int_{\h_{\R}^{\star}} \big| \Fourier_{\h_{\R}}(f_0 A_{\delta}
(\xi)) \big|^{q'} \prod_{\alpha \in \Root^+} |\langle \alpha, \xi
\rangle|^{2-q'} d \xi \Big)^{1/q'}}{\displaystyle \Big(
\int_{\h_{\R}} |f_0 A_{\delta} (x)|^q \,\ \big| k^{|\Root^+|}
A_{\delta}(x/k) \big|^{2-q} dx \Big)^{1/q}}$$ for $1 < q \le 2$,
where the supremum runs over the family of functions $f_0: \h_{\R}
\rightarrow \C$, supported in $\mt$ and symmetric under the
reflections generating the Weyl group of $G$. If
$\mathcal{K}_{f_0}(G,q)$ denotes the expression given above, then
one easily gets that $\mathcal{K}_{f_0}(G,q)$ equals
$$\frac{|\Weyl|^{\tau}}{(2 \pi)^{\tau |\Root^+|}
\mathrm{V}_G^{1/q'}} \prod_{\alpha \in \Root^+} |\langle \alpha,
\delta \rangle|^{\tau} \frac{\displaystyle \Big(
\int_{\h_{\R}^{\star}} \big| \Fourier_{\h_{\R}}(f_0 A_{\delta}
(\xi)) \big|^{q'} \prod_{\alpha \in \Root^+} |\langle \alpha, \xi
\rangle|^{2-q'} d \xi \Big)^{1/q'}}{\displaystyle \Big(
\int_{\h_{\R}} |f_0 A_{\delta} (x)|^q \prod_{\alpha \in \Root^+}
|\langle \alpha, x \rangle|^{2-q} dx \Big)^{1/q}}.$$ Moreover,
taking $q = 2$ and by Plancherel theorem on compact groups, it
follows that $\mathrm{V}_G = 1$. The boundedness of this
expression can be regarded as a weighted Hausdorff-Young
inequality of Pitt type, see \cite{B2} for more on this topic.
\end{remark}

As we pointed out in the introduction, the growth of
$\Constant_q^1(l^{p'}(n),G)$ remains open for $1 \le p < q \le 2$.
We end this paper with some remarks about this problem.

\begin{remark}
In theorem \ref{p<2} we found an extremal function $\Phi_n =
(\varphi_1, \varphi_2, \ldots \varphi_n)$, such that
$$\Constant_q^1(l^p(n),G) \ge
\frac{\|\widehat{\Phi}_n\|_{\Lebesgue_{l^p(n)}^{q'}(\G)}}
{\|\Phi_n\|_{L_{l^p(n)}^q(G)}} \ge \mathcal{K}(G,q) \,\ n^{1/p -
1/q}.$$ Our functions $\varphi_1, \varphi_2, \ldots \varphi_n$
satisfied two crucial properties, namely
\begin{itemize}
\item[$(\mathrm{P}1)$] The norm of
$\widehat{\varphi}_k(\pi)$ on $S_{d_{\pi}}^{q'}$ does not depend
on $k$ for any $\pi \in \G$.
\item[$(\mathrm{P}2)$] $\varphi_1,
\varphi_2, \ldots, \varphi_n$ have pairwise disjoint supports on
$G$.
\end{itemize}
The idea was to compare the norms of $\widehat{\Phi}_n$ and
$\Phi_n$ with $n^{1/p}$ and $n^{1/q}$ respectively. To this end,
properties $(\mathrm{P}1)$ and $(\mathrm{P}2)$ were the conditions
to be required since they provided suitable simplifications for
the original expressions of such norms. Now, if we replace
$l^p(n)$ by $l^{p'}(n)$ in the relation above, we want to compare
the norms of $\widehat{\Phi}_n$ and $\Phi_n$ with $n^{1/q'}$ and
$n^{1/p'}$ respectively. Notice that $1/p - 1/q = 1/q' - 1/p'$.
For that, we require these other properties on $\varphi_1,
\varphi_2, \ldots \varphi_n$
\begin{itemize}
\item[$(\mathrm{P}3)$] The absolute value
$|\varphi_k(g)|$ does not depend on $k$ for any $g \in G$.
\item[$(\mathrm{P}4)$] $\widehat{\varphi}_1,
\widehat{\varphi}_2, \ldots, \widehat{\varphi}_n$ have pairwise
disjoint supports on $\G$.
\end{itemize}
In the introduction we recalled that the growth of
$\Constant_q^1(l^p(n),G)$ and $\Constant_q^1(l^{p'}(n),G)$ can be
understood as dual problems with respect to the Fourier transform
operator. Now, these properties justify this point. Assuming
properties $(\mathrm{P}3)$ and $(\mathrm{P}4)$, we get
$\Constant_q^1(l^{p'}(n),G) \ge \mathcal{K}'(G,q,n) \,\ n^{1/q' -
1/p'}$, where $\mathcal{K}'(G,q,n)$ is given by
$$\mathcal{K}'(G,q,n) = \left( \frac{1}{n} \,\ \sum_{k=1}^n \left[
\frac{\|\widehat{\varphi}_k\|_{\Lebesgue^{q'}(\G)}}
{\|\varphi_k\|_{L^q(G)}} \right]^{q'} \right)^{1/q'}.$$ Hence, if
we define $\mathcal{K}'(G,q) = \inf_{n \ge 1}
\mathcal{K}'(G,q,n)$, it remains to see that $\mathcal{K}'(G,q) >
0$. We do not know if this inequality holds for any compact
semisimple Lie group and any $1 \le q \le 2$.
\end{remark}

\begin{remark}
We do not know if properties $(\mathrm{P}3)$ and $(\mathrm{P}4)$
are compatible. However, given $f_0 \in L^2(G)$ continuous and any
sequence of positive numbers $\{\varepsilon_n: \,\ n \ge 1\}$
decreasing to $0$, it is not difficult to see that there exists a
system $\Phi = \{ \varphi_n: \,\ n \ge 1\}$ of trigonometric
polynomials on $G$ satisfying
\begin{enumerate}
\item The functions $\widehat{\varphi}_1, \widehat{\varphi}_2,
\ldots$ have pairwise disjoint supports on $\G$.
\item The estimate $|\varphi_n| \le |f_0| + \varepsilon_n$ holds
in $G$.
\item The estimate $|\varphi_n| \ge |f_0| - \varepsilon_n$ holds
outside $\Omega_n$, where $\mu(\Omega_n) \rightarrow 0$ as $n
\rightarrow \infty$.
\end{enumerate}
\end{remark}

\begin{remark}
As it is well-known, $\Constant_q^1(l^{p'}(n),G) = n^{1/p-1/q}$
for any compact abelian group $G$. This equality follows by taking
$\varphi_1, \varphi_2, \ldots \varphi_n$ to be a collection of $n$
pairwise distinct characters. This motivates us to see what
happens when we consider the irreducible characters of a compact
semisimple Lie group. Let $\chi_{\lambda}$ be the character of the
irreducible representation $\pi_{\lambda}$, let us consider the
function $\Phi_n(g) = (d_{\lambda_1}^{\tau}
\chi_{\lambda_1}(g),d_{\lambda_2}^{\tau} \chi_{\lambda_2}(g),
\ldots d_{\lambda_n}^{\tau} \chi_{\lambda_n}(g))$, where
$\lambda_1, \lambda_2, \ldots \lambda_n$ are pairwise distinct
dominant weights and $\tau = 1 - 2/q'$. Then we have
$$\|\widehat{\Phi}_n\|_{\Lebesgue_{l^{p'}(n)}^{q'}(\G)} = \Big(
\sum_{k=1}^n d_{\lambda_k} \|d_{\lambda_k}^{\tau}
\widehat{\chi}_{\lambda_k}(\pi_k)\|_{S_{d_{\lambda_k}}^{q'}}^{q'}
\Big)^{1/q'} = n^{1/q'}.$$ On the other hand, applying
consecutively the Weyl integration formula and the Weyl character
formula, we get $$\|\Phi_n\|_{L_{l^{p'}(n)}^{q}(G)} = \Big(
\frac{1}{|\Weyl|} \int_{\MT} \Big( \sum_{k=1}^n |d_{\lambda_k}
A_{\delta}(t)|^{\tau p'} |A_{\lambda_k + \delta}(t)|^{p'}
\Big)^{q/p'} dm(t) \Big)^{1/q}.$$ However these relations do not
provide optimal growth. For instance, in the simplest case $G =
SU(2)$ it can be checked that there exists a constant
$\mathcal{K}_{p,q}$ depending on $p$ and $q$ such that
$$\frac{1}{n^{1/p'}} \Big( \int_{SU(2)}
\|\Phi_n(g)\|_{l^{p'}(n)}^q d \mu(g) \Big)^{1/q} \ge
\mathcal{K}_{p,q} \,\ n^{\tau}.$$
\end{remark}

\begin{remark}
If we try to find out why our attempts to get optimal growth have
failed, we need to revisit the proof of theorem \ref{p<2}. The
point is that we required the functions $\varphi_1, \varphi_2,
\ldots \varphi_n$, not only to satisfy properties $(\mathrm{P}1)$
and $(\mathrm{P}2)$, but also to be translations of a common
function. This was essential in section \ref{Section-p<2} and here
the obstacle lies in the fact that we can not take translations
since the dual object has not a group structure. This is the main
difference with the abelian case where, since the dual object is a
group, multiplication by a character in $G$ becomes a translation
in the other side of the Fourier transform operator.
\end{remark}

\begin{remark}
The quantized Rademacher system associated to a probability space
$(\Omega, \mathcal{M}, \mu)$, an index set $\Sigma$, and a family
$\{d_{\sigma}: \sigma \in \Sigma\}$ of positive integers is
defined by a collection $\mathcal{R} = \{\rho^{\sigma}: \Omega
\rightarrow O(d_{\sigma})\}_{\sigma \in \Sigma}$ of independent
random orthogonal matrices, uniformly distributed on the
orthogonal group $O(d_{\sigma})$. In \cite{GP2} we define the
notions of $\mathcal{R}$-type, $\mathcal{R}$-cotype and strong
$\mathcal{R}$-cotype of an operator space $E$. Moreover, we show
that $$\begin{array}{lcl} \mbox{Fourier type $p$} & \Rightarrow &
\mbox{strong} \ \ \mbox{$\mathcal{R}$-cotype $p'$}
\\ \mbox{Fourier cotype $p'$} & \Rightarrow &
\mbox{$\mathcal{R}$-type $p$.} \end{array}$$ This implications
allow us to work with the quantized Rademacher system where other
techniques are available to study the growth of
$\Constant_q^1(l^{p'}(n),G)$.
\end{remark}

\begin{remark}
Of course, the growth of $\Constant_q^1(l^{p'}(n),G)$ is trivially
optimal when we work with compact groups with infinitely many
inequivalent irreducible representations of the same degree $d_0$.
The unitary groups $U(n)$ are the simplest non-commutative
examples of this degenerate case. Also, it is not difficult to
check that $\Constant_2^1(l^{p'}(n),G) = n^{1/2 - 1/p'}$ by the
Plancherel theorem for compact groups.
\end{remark}

\bibliographystyle{amsplain}

\begin{thebibliography}{10}
\bibitem {An} M.E. Andersson, \emph{The Hausdorff-Young
inequality and Fourier type}, Ph. D. Thesis, Uppsala ($1993$).
\bibitem {Ba} K.I. Babenko, \emph{An inequality in the theory of
Fourier integrals}, Izv. Akad. Nauk $SSSR$ \textbf{25} ($1961$),
$531-542$.
\bibitem {B1} W. Beckner, \emph{Inequalities in Fourier analysis},
Ann. of Math. $(2)$ \textbf{102} ($1975$), $159-182$.
\bibitem {B2} W. Beckner, \emph{Pitt's inequality and the
uncertainty principle}, Proc. Amer. Math. Soc. \textbf{123}
($1995$), no. $6$, $1897-1905$.
\bibitem {ER2} E.G. Effros and Z.J. Ruan, \emph{Operator spaces},
London Math. Soc. Monogr. \textbf{23}, Oxford Univ. Press
($2000$).
\bibitem {Fo} G.B. Folland, \emph{A Course in Abstract Harmonic
Analysis}, Stud. Adv. Math., CRC Press ($1995$).
\bibitem {FH} W. Fulton and J. Harris, \emph{Representation
Theory: A First Course}, Grad. Texts in Math., Springer-Verlag,
$1991$.
\bibitem {GP} J. Garc\'{\i}a-Cuerva, J. Parcet, \emph{Vector-valued
Hausdorff-Young inequality on compact groups}. To appear in Proc.
London Math. Soc.
\bibitem{GP2} J. Garc\'{\i}a-Cuerva, J. Parcet, \emph{Quantized
orthonormal systems: A non-commutative Kwapie\'n theorem}, Studia
Math. \textbf{155} (2003), 273-294.
\bibitem {Ku} R.A. Kunze, \emph{$L_p$ Fourier transforms on
locally compact unimodular groups}, Trans. Amer. Math. Soc.
\textbf{89} ($1958$), $519-540$.
\bibitem {P1} G. Pisier, \emph{The Operator Hilbert Space OH,
Complex Interpolation and Tensor Norms}, Mem. Amer. Math. Soc.
\textbf{122} ($1996$), $1-103$.
\bibitem {P2} G. Pisier, \emph{Non-commutative vector valued
$L_p$-spaces and completely $p$-summing maps}, Ast\'{e}risque (Soc.
Math. France) \textbf{247} ($1998$), $1-111$.
\bibitem {Ru} B. Russo, \emph{The norm of the $L\sp{p}$-Fourier
transform on unimodular groups}, Trans. Amer. Math. Soc.
\textbf{192} ($1974$), $293-305$.
\bibitem {Si} B. Simon, \emph{Representations of Finite and
Compact Groups}, Grad. Stud. Math. \textbf{10}, Amer. Math. Soc.
($1996$).
\bibitem {Sj} P. Sj\"{o}lin, \emph{A remark on the Hausdorff-Young
inequality}, Proc. Amer. Math. Soc. \textbf{123} ($1995$),
$3085-3088$.

\end{thebibliography}

\end{document}